\documentclass[11pt]{article}

\usepackage{amsmath}   
\usepackage{amssymb}

\newcommand{\ed}{\end{document}}

\newcommand{\calC}{{\cal C}}
\newcommand{\calF}{{\cal F}}
\newcommand{\calR}{{\cal R}}
\newcommand{\Cn}{\calC_n}
\newcommand{\oa}{\overline{a}}
\newcommand{\ob}{\overline{b}}
\newcommand{\oab}{\overline{ab}}
\newcommand{\diag}{\mbox{\rm diag}}

\begin{document}

\renewcommand{\arraystretch}{1.2}  
\setlength{\arraycolsep}{.05in}    

\title{Universal Similarity
Factorization Equalities over 
Complex Clifford Algebras}
\author{Yongge Tian \\
Department of Mathematics and Statistics \\
Queen's University \\
Kingston,  Ontario,  Canada K7L 3N6\\
{\tt e-mail:ytian@mast.queensu.ca}}
\date{}
\maketitle


\begin{abstract}
A set of valuable universal similarity factorization equalities are established over complex Clifford algebras $\Cn.$ Through
them matrix representations of complex Clifford algebras 
$\Cn$ can directly be derived, and their properties can easily be determined.\\

\noindent {\em AMS Subject Classification}: 15A23, 15A66. \\
\noindent {\em Keywords: } Clifford algebras, matrix representations,
universal similarity factorizations. \par
\end{abstract}



\section{Introduction}
Let $\Cn$ be the complex Clifford algebra, 
\index{Clifford algebra!complex}
with the identity 1, defined on $n$ generators 
$e_1, e_2, \cdots, e_n$ subject to the multiplication laws:
\begin{align}
e_i^2 &= -1,  \quad i = 1,2,\cdots,n,                      \label{eq:eq1}\\ 
e_ie_j + e_je_i &= 0, \quad i \neq j,\;i,j=1,2,\cdots, n,  \label{eq:eq2} 
\end{align}
and $e_1e_2 \cdots e_n \neq \pm 1.$ In that case $\Cn$ is spanned as a $2^n$-dimensional vector space with $2^n$ basis
$\{e_A\},$ where the multi index $A$ ranges all naturally ordered subsets of the first positive integer set $\{1,2,\cdots,n\};$
the basis element $e_A,$ where 
$A=(i_1,i_2,\cdots,i_k)$ with $1 \leq i_1 < i_2 < \cdots < i_k \leq n,$ is defined as the product
$$
e_A \equiv e_{(i_1, i_2, \cdots, i_k)} = e_{i_1}e_{i_2} \cdots e_{i_k}, \quad 
e_{A = \emptyset} = 1. 
$$ 
For simplicity, the volume element $e_{12\cdots n}=e_1e_2 \cdots e_n$ of $\Cn$ will be denoted by $e_{[n]}$ in the
sequel. The square of the volume element is
$$
e_{[n]}^2 = (-1)^{\frac12 n(n-1)}e_1^2e_2^2\cdots e_n^2.
$$
In that case, all $a \in \Cn$ can be expressed as
\index{Clifford element!general} 
$$
a = \sum_A a_Ae_A, \quad  a_A \in \calC, 
$$   
where $A$ ranges all naturally ordered subsets of $\{1,2,\cdots,n\}.$ We shall adopt the  following notation from now on:
$\Cn := { \cal C }\{ \,  e_1,\cdots, \ e_n  \}.$  

Clifford algebras (real or complex) have been studied for many years and their algebraic properties are well known. In particular, all Clifford algebras are classified as matrix algebras, or as direct sums of matrix algebras  over the fields of real or complex numbers, or the quaternion ring (see, e.g., \cite{At,1,Ka,O1,O2,O3,4,5}). For complex Clifford algebras $\Cn,$ it is well known that they can faithfully be realized as 
\index{generators}%
\index{multiplication laws}%
\index{volume element}%
\index{basis element}%
\index{Clifford algebra!complex}%
\index{matrix algebra!total complex}%
certain matrix algebras over ${\cal C},$ and the general algebraic isomorphism is
\begin{equation} 
\Cn \simeq 
\begin{cases}
  \calC(2^{\frac{n}{2}})        & \text{if $n$ is even,}  \\[1ex]
^2{\cal C}(2^{\frac{n-1}{2}})   & \text{if $n$ is odd,}
\end{cases}
\label{eq:eq3}
\end{equation}
where $\calC(s) $ stands for the $ s \times s $ total complex matrix 
algebra, and $ ^2{\cal C}(s) $ stands for the complex matrix algebra
$$ 
  ^2\calC(s) = \left\{ \left. \left[\begin{array}{cc}  
                                     A & 0 \\ 
                                     0 & B 
                                    \end{array} 
                              \right] \
               \right|  \ A, B \in \calC(s) \right\}. 
$$
\index{algebraic isomorphism}%
In this article we improve this relationship to a new level by establishing a set of valuable universal similarity  factorization
equalities between elements of $\Cn$ and complex 
matrices over $ \calC(2^{\frac{n}{2}})$ or $\calC(2^{\frac{n-1}{2}})$ for all $n.$ Through these universal factorization
equalities the complex matrix representations of elements in $\Cn$ can explicitly be established. Based on them various
results in the complex matrix theory can directly  be extended to complex Clifford numbers.
\index{universal similarity factorization equality}%
\index{matrices!complex}%
\index{Clifford number!complex}%
Moreover, one can also easily develop matrix analysis over  complex Clifford algebras.

We first establish some basic universal similarity factorization equalities for elements in Clifford algebras with dimensions
$2$ and $4,$ which will be directly applied to establish  some more general results for $\Cn.$

\medskip
\noindent {\bf Lemma 1.} Let $\calF$ be an algebraically closed field and 
$\calF_1 = \calF\{e \ | \ e^2 = u \}$ be a Clifford algebra defined on a generator $e$ with 
$e^2 = u \in \calF$ and $u \neq 0.$ Then all $a \in \calF_1$ can be written as 
$a = a_0 + a_1e,$ where $a_0,a_1 \in {\cal F}.$ Moreover, define 
$\oa = a_0 - a_1e.$ In that case, $a$ and $\oa$ satisfy the following universal
similarity factorization equality
\begin{equation}
P \left[\begin{array}{rr} a & 0 \\ 
                          0 & \oa 
        \end{array}
  \right]
P^{-1}
= \left[\begin{array}{cc}
        a_0 + \sqrt{u}a_1 & 0  \\ 
                        0 & a_0 - \sqrt{u}a_1 
        \end{array} \right], 
\label{eq:eq4}
\end{equation} 
where $P$ and $P^{-1}$ have the following universal forms (no relation with~$a):$
\begin{align} 
P &= \frac12 \left[\begin{array}{cc}
                   1 + \frac{1}{\sqrt{u}}e   & -( \sqrt{u} - e ) \\ 
                   \frac{1}{u}(\sqrt{u} - e) & 1+ \frac{1}{\sqrt{u}}e
                  \end{array} 
            \right],
\notag \\
P^{-1} &= \frac12 \left[\begin{array}{cc} 
                        1 + \frac{1}{\sqrt{u}}e    & \sqrt{u} - e  \\
                        -\frac{1}{u}(\sqrt{u} - e) & 1+ \frac{1}{\sqrt{u}}e
                       \end{array}
                 \right]. 
\label{eq:eq5}
\end{align}

\noindent {\bf Proof.} Note that $e^2=u.$ It is easy to verify that 
$$
\left[ \begin{array}{cc}  
        1 & e  \\ 
        e^{-1} & -1 
       \end{array} 
\right]
\left[ \begin{array}{cr}
        a & 0 \\
        0 & \oa \end{array} \right]
\left[ \begin{array}{cc}  
        1  & e \\ 
        e^{-1} & -1 
       \end{array} 
\right]
= 
\left[ \begin{array}{cc}  
        a_0 & ua_1 \\ 
        a_1 & a_0 \oa
       \end{array} 
\right].
$$
On the other hand, it is also easy to verify that
$$
\left[ \begin{array}{cc} 
        1 & \sqrt{u} \\ 
        \frac{1}{\sqrt{u}} & -1 
       \end{array} 
\right] 
\left[ \begin{array}{cc} 
        a_0 + ua_1 & 0 \\
        0 & a_0 - ua_1 
       \end{array} 
\right]
\left[ \begin{array}{cc} 
        1 & \sqrt{u} \\ 
        \frac{1}{\sqrt{u}} & -1 
       \end{array} 
\right]
= 
\left[ \begin{array}{cc} 
        a_0 & ua_1 \\ 
        a_1 & a_0 \oa 
       \end{array} 
\right].
$$
Combining the above two equalities yields (\ref{eq:eq4}) and (\ref{eq:eq5}). \quad $\square$

\medskip
\noindent{\bf Lemma 2.} Let $M_2(\calF)$ be the $2 \times 2$ total matrix algebra over an arbitrary field $\calF$ with its
basis satisfying the multiplication rules
\begin{equation}
\tau_{pq} \tau_{st} = 
\begin{cases}
\tau_{pt}, & \text{$q = s,$} \\
0,         & \text{$q \neq s,$}
\end{cases}
\label{eq:eq6}
\end{equation}
for $p,q,s,t = 1,2.$ Then all 
$a=a_{11}\tau_{11}+a_{12}\tau_{12}+a_{21}\tau_{21}+a_{22}\tau_{22} \in M_2(\calF),$ where $a_{pq} \in \calF,$
satisfy the following universal similarity factorization equality
\begin{equation}
Q \left[ \begin{array}{cc} 
          a & 0 \\
          0 & a  
         \end{array} 
  \right] 
Q^{-1}
= 
\left[ \begin{array}{cc}
        a_{11} & a_{12} \\  
        a_{21} & a_{22}  
       \end{array} 
\right],
\label{eq:eq7}
\end{equation} 
where $Q$ has the following universal form  
\begin{equation} 
Q =  Q^{-1} = \left[ \begin{array}{cccc} 
                      \tau_{11} & \tau_{21} \\
                      \tau_{12} & \tau_{22}  
              \end{array} 
\right].
\label{eq:eq8}
\end{equation}  

\medskip
\noindent {\bf Proof.}  Follows directly from a verification. \quad $\square$  

\medskip
\noindent {\bf Lemma 3.} Let $\calF$ be an algebraically closed field of characteristic not two, and $\calF_2 = \calF\{e_1,
e_2\}$ be a Clifford algebra defined on two generators $e_1,e_2$ with $e^2_1 = u \in \calF, e^2_2 = v \in \calF$ and $u
\neq 0, v \neq 0.$ Then all $a \in \calF_2$ can be written as
\begin{equation}
a = a_0 + a_1e_1 + a_2e_2 + a_3e_{12}, 
\label{eq:eq9}
\end{equation} 
where $a_0,\ldots,a_3 \in \calF.$ In that case, $aI_2 $ satisfies the following universal similarity 
factorization  equalities
\begin{equation}
R \left[ \begin{array}{rr}  
          a & 0 \\ 
          0 & a
         \end{array} 
\right] 
R^{-1}
= 
\left[ \begin{array}{cc} 
        a_0 + \sqrt{u}a_1 & v( a_2 + \sqrt{u}a_3 ) \\
        a_2 - \sqrt{u}a_3 & a_0 - \sqrt{u}a_1 
       \end{array} 
\right] ,
\label{eq:eq10}
\end{equation} 
and
\begin{equation}
 T \left[ \begin{array}{rr} 
           a & 0 \\ 
           0 & a 
          \end{array}
   \right] 
T^{-1}
= 
\left[ \begin{array}{cc} 
        a_0 + \sqrt{v}a_1 & u( a_1 - \sqrt{v}a_3 ) \\  
        a_1 + \sqrt{v}a_3 & a_0 - \sqrt{v}a_1 
       \end{array} 
\right],
\label{eq:eq11}
\end{equation} 
where
\begin{equation} 
R = R^{-1} = \frac12 \left[ \begin{array}{cc} 
    1 + \frac{1}{\sqrt{u}}e_1 &  e_2 - \frac{1}{\sqrt{u}}e_{12}  \\                                 
 \frac{1}{v}( e_2 + \frac{1}{\sqrt{u}}e_{12} ) & 1 - \frac{1}{\sqrt{u}}e_1
                            \end{array} \right],
\label{eq:eq12}
\end{equation} 
and
\begin{equation}
 T = T^{-1}  = \frac12 \left[ \begin{array}{cc} 
     1 + \frac{1}{\sqrt{v}}e_2 & e_1  + \frac{1}{\sqrt{v}}e_{12}  \\
     \frac{1}{u}( e_1 - \frac{1}{\sqrt{v}}e_{12}) & 1 - \frac{1}{\sqrt{v}}e_2  
                              \end{array} \right].
\label{eq:eq13}
\end{equation} 

\medskip
\noindent {\bf Proof.} By Lemma 2, we take the change of basis of $\calF_2$ as follows
\begin{subequations}\label{eq:eq14}
\begin{align} 
\tau_{11} &= \frac12\left(1+ \frac{1}{\sqrt{u}}e_1\right),
& \tau_{12} &= \frac{1}{2v}\left(e_2 + \frac{1}{\sqrt{u}}e_{12} \right), \label{eq:14a} \\[1ex]
\tau_{21} &= \frac12 \left(e_2 - \frac{1}{\sqrt{u}} e_{12}\right), 
& \tau_{22} &= \frac12 \left(1 -\frac{1}{\sqrt{u}} e_1 \right).          \label{eq:14b} 
\end{align}
\end{subequations}
Then it is not hard to verify that this new basis satisfies the multiplication rules in (\ref{eq:eq6}). In this new 
basis, any 
$a = a_0 + a_1e_1 + a_2e_2 + a_3e_{12} \in \calF_2$ can be rewritten as
$$
a=(a_0 + \sqrt{u} a_1)\tau_{11} + (va_2 +\sqrt{u} a_3)\tau_{12} +
  (a_2 - \sqrt{u} a_3)\tau_{21} + (a_0 - \sqrt{u}a_1)\tau_{22}.
$$
Substituting this and (\ref{eq:eq14}) into (\ref{eq:eq7}) and
(\ref{eq:eq8}), we can obtain 
(\ref{eq:eq10}) and (\ref{eq:eq12}). By the similar approach, we can get (\ref{eq:eq11}) and 
(\ref{eq:eq13}). \quad $ \square $  

\section{Main results}
Notice that $\calC$ is an algebraically closed field. We then can
establish a set of universal similarity factorization equalities for
elements in $\Cn$ on the basis of Lemmas 1 and 3. 

\medskip
\noindent {\bf Theorem 4.}  Let $a \in \calC_1 = \calC\{e\}$ be given. Then $a$ can be written as $a=a_0 + a_1e,$ where
$a_0,a_1 \in \calC = \{x + iy \, | \, x,y \in \calR\}.$  Moreover denote 
by $\oa = a_0 - a_1e,$ called the conjugate of $a.$ 
\index{conjugate}%
Then $a$ and $\oa $ satisfy the following universal similarity factorization equality 
\begin{equation}
P_1 \left[ \begin{array}{rr} 
            a & 0 \\ 
            0 & \oa  
           \end{array} 
    \right] 
P_1^{-1}
= 
\left[ \begin{array}{cc} 
        a_0 + a_1i & 0  \\ 
        0 & a_0 - a_1i 
       \end{array} 
\right],
\label{eq:eq15}
\end{equation}
where $P_1$ and $P_1^{-1}$ have the following universal forms (no relation with~$a):$
\begin{equation}
P_1 = \frac12 \left[ \begin{array}{cc} 
                     1 - ie & -(i - e) \\ 
                     -(i - e) & 1 - ie  
                     \end{array} 
              \right],
\quad
P_1^{-1} = \frac12 \left[ \begin{array}{cc}
                           1 - ie &  i - e \\ 
                           i - e  &  1 - ie
                          \end{array} \right]. 
\label{eq:eq16}
\end{equation}

\medskip 
\noindent {\bf Proof.} Let $\calF = \calC$ and $u = -1$ as in Lemma 1. Then 
(\ref{eq:eq15}) and (\ref{eq:eq16}) follow directly from (\ref{eq:eq4}) and (\ref{eq:eq5}). 
\quad $\square$  

\medskip
It is easy to verify that 
\begin{equation}
\overline{\oa} = a,  \quad  \overline{a + b} = \oa + \ob, \quad
\oab=\oa\ob, \quad \overline {\lambda a} = \overline { a \lambda } =
 \lambda \oa. 
\label{eq:eq15a}
\end{equation}  
hold for all $a,b \in \calC_1$ and $\lambda \in \calC.$ According to 
(\ref{eq:eq15}), we define a map from $\calC_1$  to the double field  $ \calC
 \oplus \calC $ by  
\begin{equation}
\phi_1: a  = a_0 + a_1e \in \calC_1 \longrightarrow \left[ \begin{array}{cc} 
        a_0 + a_1i & 0  \\ 
        0 & a_0 - a_1i 
       \end{array} 
\right] \in  \calC \oplus \calC.
\label{eq:eq15b}
\end{equation}
Then it is easy to derive from (\ref{eq:eq15}) and  (\ref{eq:eq15a}) the following 
properties.  

\medskip
\noindent {\bf Corollary 5.} Let $a=a_0+a_1e,b= b_0 + b_1e \in \calC_1, 
\lambda \in \calC$ be given, and $\phi_1$ be defined by (\ref{eq:eq15b}).  Then 
\begin{itemize}
\item[(a)] $a=b \Longleftrightarrow \phi_1(a) = \phi_1(b).$ 
\item[(b)] $\phi_1(a + b ) =  \phi_1(a) + \phi_1(b),\,
            \phi_1(ab)=\phi_1(a)\phi_1(b),\, 
            \phi_1(\lambda a) = \lambda \phi_1(a),$ 
           $\phi_1(1)= I_2.$
\item[(c)] $\phi_1(\oa) = \left[ \begin{array}{cc} 
                                  0 &  1 \\ 
                                  1 &  0  
                                 \end{array} 
                          \right]
\phi_1(a) \left[ \begin{array}{cc} 
                  0 & 1 \\ 
                  1 &  0  
                 \end{array}
\right].$
\item[(d)] Denote $a^{\#} = \overline{a_0} - \overline{a_1}e,$ then 
           $\phi_1(a^{\#}) = \phi_1^*(a),$ the conjugate transpose of the complex matrix $\phi_1(a).$ 
\index{conjugate transpose}%
\item[(e)] $a = \frac{1}{4}[1 - ie,i - e] \phi_1(a)[1 -ie,e - i]^T.$
\item[(f)] $\det [\phi_1(a)] =  a_0^2 + a_1^2.$ 
\item[(g)] $a$ is invertible $ \Longleftrightarrow $ $\phi_1(a)$ is invertible, in which
case $\phi_1(a^{-1}) = \phi_1^{-1}(a).$ 
\end{itemize}
The properties in Corollary 5(a) and (b) show that through the bijective map (\ref{eq:eq15b}),
\index{bijective map}%
the Clifford algebra $\calC_1$ is algebraically isomorphic to the  double field  
$ \calC \oplus \calC $, and  $\phi_1(a)$ is a faithful matrix  representation of $a$ 
in $ \calC \oplus \calC.$

Notice that $P_1$ and $P_1^{-1}$ in (\ref{eq:eq15}) have no relation with $a.$ Thus the equality in Theorem 4 
can also be extended to all matrices over the complex Clifford algebra $\calC_1.$ 

\medskip
\noindent {\bf Theorem 6.} Let $A = A_0 + A_1e \in \calC_1^{m \times n}$ be given, where  
$A_0, A_1 \in \calC^{m \times n}.$ Then $A$ and its conjugate $\overline{A} = A_0 - A_1e$ satisfy the following
universal factorization equality
\begin{align}
J_{2m} \left[ \begin{array}{cc}  
               A & 0 \\ 
               0 & \overline{A}
              \end{array} 
       \right] 
J_{2n}^{-1}
&= 
\left[ \begin{array}{cc} 
        A_0 + A_1i & 0 \\ 
        0 & A_0 - A_1i
       \end{array} \right],
\label{eq:eq17}
\end{align}
where  
\begin{align} 
J_{2m} &= \frac12 \left[ \begin{array}{cc} 
                         (1 -ie)I_m   & -(i - e )I_m \\ 
                         -(i - e)I_m & (1 - ie_1)I_m
                        \end{array} 
                 \right],\notag \\
J_{2n}^{-1} &= \frac12 \left[ \begin{array}{cc} 
                              (1 - ie_1)I_n  & ( i - e )I_n \\ 
                             ( i - e )I_n & (1- ie)I_n  
                             \end{array} 
                      \right].
\label{eq:eq18}
\end{align}
In particular, if $m = n,$ then (\ref{eq:eq17}) becomes a universal similarity factorization  equality over $\calC_1$. 

\medskip
It is easy to verify that for any matrices  $A,B \in \calC_1^{m \times n }, C \in \calC_1^{n \times p },$ 
and $\lambda \in \calC$
\begin{equation}
\overline{\overline{A}} = A,  \quad  \overline{A + B} = \overline{A} + \overline{B}, \quad
 \overline{AC}=\overline{A}\,\overline{C}, \quad \overline{\lambda A} = \overline{ A \lambda } =
 \lambda \overline{A}. 
\label{eq:eq18a}
\end{equation}  
Now according to (\ref{eq:eq17}), we define the complex matrix representation of a matrix $ A = A_0 + A_1e \in \calC_1^{m \times n}$ by  $\Phi_1(A) = \left[\begin{array}{cc} A_0 + A_1i & 0 \\  0 & A_0 - A_1i \end{array} \right].$ Then the following properties can  be easily derived from (\ref{eq:eq17}) and (\ref{eq:eq18a}).

\medskip
\noindent {\bf Corollary 7.} Let $A,B \in \calC_1^{m \times n}, C \in \calC_1^{n \times p}$
and $\lambda \in \calC$ be given. Then
\begin{itemize}

\item[(a)] $A=B \Longleftrightarrow \Phi_1(A) = \Phi_1(B).$ 

\item[(b)] $\Phi_1(A + B ) = \Phi_1(A) + \Phi_1(B).$ 

\item[(c)] $\Phi_1(AC) = \Phi_1(A) \Phi_1(C),\;
            \Phi_1(\lambda A) = \lambda \Phi_1(A), \;
            \Phi_1(I_m )= I_{2m}. $

\item[(d)] $\Phi_1(\overline{A}) = \left[ \begin{array}{cc} 
                                          0 & I_m \\
                                          I_m & 0  
                                         \end{array} 
                                  \right]
\Phi_1(A) 
\left[ \begin{array}{cc} 
        0 & I_n \\ 
        I_n & 0  
       \end{array}
\right]. $

\item[(e)] Let $A = A_0 + A_1e$ and denote $A^{\#} =  A_0^* - A_1^*e,$ where $ A_0^*$ and $ A_1^*$ are the conjugate
transposes of the complex matrices $ A_0$ and $A_1.$ Then 
$\Phi_1(A^{\#}) = \Phi_1^*(A),$ the conjugate transpose of the complex matrix $\Phi_1(A). $

\item[(f)] $A = \frac{1}{4}[(1 - ie)I_m,(i - e )I_m] \Phi_1(A)[(1 -ie)I_n,(e - i)I_n]^T.$

\item[(g)] $A$ is invertible $\Longleftrightarrow $ $\Phi_1(A)$ is invertible, in which case 
$\Phi_1(A^{-1}) = \Phi_1^{-1}(A).$  

\item[(h)] $p_A(A)=0,$ where $p_A(x) = \det [x I_{2m} - \Phi_1(A)].$ 
\end{itemize}

\medskip
\noindent {\bf Theorem 8.}  Let $a \in \calC_{2} = \calC \{e_1,e_2\},$ the  
complex quaternion algebra.
\index{quaternion algebra!complex}%
 Then $a$ can be written as
$$
a = a_0 + a_1e_1 + a_2e_2 + a_3e_{12},
$$
where $a_0,\ldots,a_3 \in \calC.$ In that case, $aI_2$ satisfies the following universal similarity factorization equality
\begin{equation}
P_{2} \left[ \begin{array}{cc} 
              a & 0 \\ 
              0 & a 
             \end{array} 
      \right]
P_{2}^{-1}
= 
\left[ \begin{array}{cc} 
        a_0 + a_1i & -(a_2 + a_3i) \\
        a_2 - a_3i & a_0 - a_1i 
       \end{array} 
\right],
\label{eq:eq19}
\end{equation}
where $P_{2}$ has the universal form  
\begin{equation}
P_{2} = P_{2}^{-1} = \frac12 \left[ \begin{array}{cc} 
                                     1 - ie_1 & e_2 + ie_{12} \\  
                                     -e_2 + ie_{12} & 1 + ie_1 
                                    \end{array} 
                             \right].  
\label{eq:eq20}
\end{equation}

\noindent {\bf Proof.}  Let $\calF = \calC$ and $u = v = -1$ in Lemma 3. Then 
(\ref{eq:eq19}) and (\ref{eq:eq20}) follow from (\ref{eq:eq10}) and (\ref{eq:eq12}). 
\quad $\square$  

\medskip
According to (\ref{eq:eq19}), define a map from $\calC_2$  to the $2 \times 2$ complex matrix algebra 
$\calC^{2 \times 2}$ by  
\begin{equation}
\phi_2: a  =  a_0 + a_1e_1 + a_2e_2 + a_3e_{12} \in \calC_2 \longrightarrow \left[ \begin{array}{cc} 
         a_0 + a_1i & -(a_2 + a_3i) \\
        a_2 - a_3i & a_0 - a_1i 
       \end{array} 
\right] \in \calC^{2 \times 2} .
\label{eq:eq20a}
\end{equation}
Then we can easily derive from (\ref{eq:eq19}) the following properties.  
 
\medskip
\noindent{\bf Corollary 9.}  Let $a, b \in \calC_2 = \calC \{e_1,e_2\}$ and $\lambda \in \calC$ be given. Then
\begin{itemize}

\item[(a)] $a=b \Longleftrightarrow \phi_2(a) = \phi_2(b).$ 

\item[(b)] $\phi_2(a + b ) =  \phi_2(a) + \phi_2(b),\; 
            \phi_2(ab) = \phi_2(a) \phi_2(b),\;
            \phi_2( \lambda a) = \lambda \phi_2(a),$ 
           $\phi_2(1 )= I_2.$

\item[(c)] Let $a=a_0 + a_1e_1 + a_2e_2 + a_3e_{12} \in \calC_2,$ and denote 
$a^{\#} = \overline{a_0} - \overline{a_1}e_1 - \overline{a_2}e_2 - \overline{a_3}e_{12}.$  Then 
$\phi_2(a^{\#}) = \phi_2^*(a),$ the conjugate transpose of the complex matrix $\phi_2(a).$

\item[(d)] $a = \frac{1}{4}[1 - ie, e_2 + ie_{12} ] \phi_2(a)[1 -ie,-e_2 + ie_{12}]^T .$

\item[(e)] $\det [ \phi_2(a) ] =  a_0^2 + a_1^2 + a_2^2 + a_3^2.$ 

\item[(f)] $a$ is invertible $\Longleftrightarrow \phi_2(a)$ is invertible, in which case
$\phi_2(a^{-1}) = \phi_2^{-1}(a).$

\item[(g)] The two elements $a$ and $b$ in $\calC_2$ are similar, i.e., there is an invertible
\index{similarity}%
$x \in \calC_2$ such that $ax = xb$ if and only if the two complex matrices $\phi_2(a )$ and $\phi_2(b)$ are similar over
$\calC.$ 
\end{itemize}

The properties in Corollary 9(a) and (b) clearly show that through the bijective map  (\ref{eq:eq20a}) 
the Clifford algebra $\calC_2,$ i.e., the complex quaternion algebra, is
algebraically isomorphic to the complex matrix algebra $\calC^{ 2 \times
2},$ and $\phi_2(a)$ is a faithful matrix representation of $a$ in 
$\calC^{ 2 \times 2}.$
\index{matrix representation!faithful}
\index{quaternion algebra!complex}

Notice that $P_2$ and $P_2^{-1}$ in (\ref{eq:eq20}) have no relation to $a.$ Thus the equality in Theorem 8 can also be
extended to all matrices over $\calC_2.$ 

\medskip
\noindent {\bf Theorem 10.}  Let 
$A = A_0 + A_1e_1 + A_2e_2 + A_3e_{12} \in \calC_{2}^{ m \times n} = 
\calC^{ m \times n } \{e_1,e_2\}$ be given where $A_0,\ldots,A_3 \in \calC^{ m \times n}.$
Then $A$  satisfies the following universal factorization equality  
\begin{align}
K_{2m} \left[ \begin{array}{cc} 
               A & 0 \\ 
               0 & A 
       \end{array}
       \right]
K_{2n}^{-1}
&= 
\left[ \begin{array}{cc} 
        A_0 + A_1i & -(A_2 + A_3i) \\
        A_2 - A_3i & A_0 - A_1i 
       \end{array} 
\right],
\label{eq:eq21}
\end{align}
where
\begin{equation}
K_{2t} = K_{2t}^{-1} = \frac12 \left[ \begin{array}{cc} 
                                       (1 - ie_1)I_t & (e_2 + ie_{12})I_t \\
                                       (-e_2 + ie_{12})I_t & (1 + ie_1)I_t 
                                      \end{array} 
                               \right].
\label{eq:eq22}
\end{equation}
In particular, if $m = n,$ then (\ref{eq:eq21}) becomes a universal similarity factorization equality over $ \calC_2.$ 

\medskip

According to (\ref{eq:eq21}), we define the complex representation of a matrix $ A = A_0 + A_1e \in \calC_2^{m \times n}$ by  $\Phi_2(A) = \left[\begin{array}{cc} A_0 + A_1i & -(A_2 + A_3i) \\ A_2 - A_3i & A_0 - A_1i  \end{array} \right].$ 
Then the following properties are easy to verify by (\ref{eq:eq21}).
 
\medskip
\noindent {\bf Corollary 11.}  Let $A,B \in \calC_2^{m \times n }, C \in \calC_2^{n \times p },$ and $\lambda \in \calC$ be
given. Then
\begin{itemize}

\item[(a)] $A = B \Longleftrightarrow \Phi_2(A) = \Phi_2(B).$ 

\item[(b)] $\Phi_2(A + B ) =  \Phi_2(A) + \Phi_2(B).$

\item[(c)] $\Phi_2(AC) = \Phi_2(A) \Phi_2(C),\; 
            \Phi_2(\lambda A) = \lambda \Phi_2(A),\;
            \Phi_2(I_m )= I_{2m}.$

\item[(d)] Let $A = A_0 + A_1e_1 + A_2e_2 +A_3e_{12}$ and denote 
$A^{\#}= A_0^* - A_1^*e_1 - A_2^*e_2 - A_3^*e_{12},$ then $\Phi_2(A^{\#}) = \Phi_2^*(A),$ the conjugate transpose of
the complex matrix $\Phi_2(A).$

\item[(e)] $A = \frac{1}{4}[(1 - ie_1)I_m,(e_2 +ie_{12})I_m] 
                \Phi_2(A)[(1 -ie_1)I_n,( -e_2 + ie_{12})I_n]^T.$

\item[(f)] $A$ is invertible $\Longleftrightarrow \Phi_2(A)$ is invertible, in which case 
$\Phi_2(A^{-1}) = \Phi_2^{-1}(A).$  

\item[(g)] $p_A(A) = 0,$ where $p_A(\lambda) = \det [\lambda I_{2m} - \Phi_2(A)],$ the characteristic polynomial of
$\Phi_2(A).$

\item[(h)] Two square matrices $A$ and $B$ are similar over ${\calC}_2,$
i.e., there is an invertible matrix $X$ over ${\cal C}_2$ such that $AX =
XB$ if and only if $\Phi_2(A )$ and $\Phi_2(B)$ are similar over $\calC.$ 
\end{itemize}

In the next several results we only present the basic universal similarity factorization equalities without listing their
operation properties and their extensions to matrices over 
$\Cn.$  

\medskip
\noindent {\bf Theorem 12.} Let $a \in \calC_{3} = \calC \{e_1,e_2,e_3\}$ be given. Then $a$ can factor as
\begin{equation}
a=a_0 + a_1e_{[3]}, 
\label{eq:eq23}
\end{equation}
where
$$
a_0, \ a_1 \in \calC_2 =\calC \{e_1,e_2\}, \quad e_{[3]}^2 = 1.
$$
Moreover, define $\oa = a_0 - a_1e_{[3]}.$ In that case, the diagonal matrix 
$D_a = \diag ( aI_2 ,\oa I_2)$ satisfies the following universal  similarity factorization equality
\begin{align}
P_{3}D_a P_{3}^{-1} &= \left[ \begin{array}{cc} 
                              \phi_{2}(a_0) + \phi_{2}(a_1) & 0  \\ 
                              0 & \phi_{2}(a_0) - \phi_{2}(a_1)
                             \end{array} \right] \notag \\[1ex]
&:= \phi_{3}(a) \in \, ^2\calC^{ 2 \times 2},
\label{eq:eq24}
\end{align}
where $\phi_{2}(a_t),\; t = 0,1,$ is the matrix representation of $a_t$ in $\calC^{2 \times 2}$ defined 
in (\ref{eq:eq20a}) and
\begin{align}
P_{3} &= \frac12 \left[ \begin{array}{cc} 
                        (1 + e_{[3]}) P_{2} & -( 1 - e_{[3]})P_{2} \\  
                        (1 - e_{[3]}) P_{2} & (1 + e_{[3]})P_{2}
                       \end{array} 
                \right], 
\label{eq:eq25} \\
P_{3}^{-1} &= \frac12 \left[ \begin{array}{cc} 
                              P_{2}^{-1}(1 + e_{[3]}) & P_{2}^{-1}( 1 - e_{[3]}) \\
                             -P_{2}^{-1}(1 - e_{[3]}) & P_{2}^{-1}(1 + e_{[3]}) 
                             \end{array} 
                      \right],
\label{eq:eq26}
\end{align}
where $P_{2}$ and $P_{2}^{-1}$ are given by (\ref{eq:eq20}). 

\medskip 
\noindent {\bf Proof.} Notice that $be_{[3]} = e_{[3]}b$ holds for all 
$b \in \calC_2 =\calC \{e_1,e_2\}.$ We have by applying (\ref{eq:eq19}) to $a$ in 
(\ref{eq:eq23}) that
$$
P_{2}(aI_2) P_{2}^{-1} =  P_{2}(a_0I_2) P_{2}^{-1} + P_{2}(a_1I_2) P_{2}^{-1} e_{[3]} = \phi_{2}(a_0) +
\phi_{2}(a_1)e_{[3]} :=  \psi(a),
$$
and
$$
P_{2}(\overline{a}I_2) P_{2}^{-1} = \phi_{2}(a_0) - \phi_{2}(a_1)e_{[3]} :=  \psi(\oa).
$$
Next we build, according to Lemma 1, a matrix and its inverse as follows 
\begin{align}
V &= \frac12 \left[ \begin{array}{cc} 
                    (1 + e_{[3]})I_2 & -( 1 - e_{[3]})I_2 \\ 
                    (1 - e_{[3]})I_2 & (1 + e_{[3]}) I_2
                   \end{array} 
            \right],
\nonumber \\
V^{-1} &= \frac12 \left[ \begin{array}{cc}
                         (1 + e_{[3]})I_2 & ( 1 - e_{[3]})I_2 \\  
                        -(1 - e_{[3]})I_2 & (1 + e_{[3]})I_2 
                        \end{array}
                 \right].
\nonumber
\end{align}
and then calculate to get 
$$
V \left[ \begin{array}{cc} 
          \psi(a) & 0 \\ 
          0 & \psi(\oa)
         \end{array} 
  \right] 
V^{-1}
=  \left[ \begin{array}{cc} 
           \phi_{2}(a_0) + \phi_{2}(a_1) & 0  \\ 
           0 & \phi_{2}(a_0) - \phi_{2}(a_1) 
          \end{array} 
   \right].
$$
Finally substituting $\psi(a) = P_{2}(aI_2) P_{2}^{-1}$ and 
$\psi(\oa) = P_{2}(\oa I_2)P_{2}^{-1}$ into the left-hand side of the above equality yields 
(\ref{eq:eq24}), (\ref{eq:eq25}), and (\ref{eq:eq26}).
\quad $\square$

\medskip
\noindent {\bf Theorem 13.}  Let $a \in \calC_{4} = \calC\{ e_1,e_2,e_3,e_4 \}$ be given. Then $a$ can factor as
\begin{equation}
a = a_0 + a_1e_{123} + a_2e_{124} + a_3e_{43} = a_0 + e_{123}a_1 + e_{124}a_2 + e_{43}a_3, 
\label{eq:eq27}
\end{equation}
where 
$$
a_0, a_1, a_2, a_3 \in \calC_{2} = \calC \{e_1,e_2\},
$$
$$ 
 e_{123}^2  = 1, \quad  e_{124}^2  = 1, \quad e_{43}= e_{123}e_{124} = - e_{124}e_{123}.
$$
In that case, $aI_4$ satisfies the following universal similarity factorization equality
\begin{align}
P_{4}( aI_4) P_{4}^{-1} &= 
\left[ \begin{array}{cc}  
        \phi_{2}(a_0) + \phi_{2}(a_1)  & \phi_{2}(a_2) + \phi_{2}(a_3) \\ 
        \phi_{2}(a_2) - \phi_{2}(a_3)  & \phi_{2}(a_0) - \phi_{2}(a_1)  
       \end{array} 
\right] \notag \\[1ex]
&:= \phi_{4}(a) \in \calC^{4\times 4},
\label{eq:eq28}
\end{align}
where $ \phi_{2}(a_t),\, t = 0, \ldots, 3,$ is the matrix representation of $a_t$ in 
$\calC^{2 \times 2}$ defined in (\ref{eq:eq20a}) and 
\begin{equation}
P_{4} = P_{4}^{-1} = \frac12 \left[ \begin{array}{cc}
                                     (1 + e_{[3]})P_{2} & ( e_{124} - e_{43}) P_{2} \\ 
                                     (e_{124} + e_{43}) P_{2} & ( 1 - e_{[3]})P_{2} 
                                    \end{array} 
                             \right],
\label{eq:eq29}
\end{equation}
where $P_{2}$ is given in (\ref{eq:eq20}). 

\medskip
\noindent {\bf Proof.}  Note that the commutative rules 
$be_{123} = e_{123}b, be_{124} = e_{124}b, be_{43} = e_{43}b$ hold for all 
$b \in \calC_{2} = \calC\{e_1,e_2\}.$ Thus, it follows from (\ref{eq:eq19})
that
\begin{align*}
&\lefteqn{P_{2}(aI_2) P_{2}^{-1}} 
\\
&= P_{2}(a_0I_2)P_{2}^{-1} + P_{2}(a_1I_2) P_{2}^{-1}e_{123} 
             + P_{2}(a_2I_2) P_{2}^{-1}e_{124} + P_{2}(a_3I_2) P_{2}^{-1}e_{43}
\\
&= \phi_{2}(a_0) + \phi_{2}(a_1)e_{123} + \phi_{2}(a_2)e_{124} + 
   \phi_{2}(a_3)e_{43} 
\\
&:= \psi(a).
\end{align*}
Next  building, according to Lemma 3, a matrix  and its inverse as follows
$$  
V = V^{-1} = \frac12 \left[ \begin{array}{cc} 
                             (1 + e_{123} )I_2 & ( e_{124} - e_{43})I_2 \\ 
                             (e_{124} + e_{43})I_2  & (1 - e_{123} )I_2
                            \end{array} 
                     \right],
$$ 
and applying them to $\psi(a)$ given above, we obtain
$$
V \left[ \begin{array}{cc}  
          \psi(a) & 0 \\ 
          0 & \psi(a)
         \end{array} 
  \right]
V^{-1}
= \left[ \begin{array}{cc}  
          \phi_{2}(a_0) + \phi_{2}(a_1) & \phi_{2}(a_2) + \phi_{2}(a_3) \\ 
          \phi_{2}(a_2) - \phi_{2}(a_3) & \phi_{2}(a_0) - \phi_{2}(a_1)
         \end{array} 
\right].
$$
Finally substituting $\psi(a) = P_{2}(aI_2)P_{2}^{-1}$ into its left-hand side yields 
(\ref{eq:eq28}) and (\ref{eq:eq29}). \quad $ \square $ 

\medskip
By induction, we have the following two  general results. 

\medskip
\noindent {\bf Theorem 14.} Suppose that there is an independent invertible matrix $P_{n}$ over 
$\calC_{n} = \calC \{e_1, \ldots, e_{n} \}$  with $n$ even such that
\begin{equation}
P_{n}( aI_{2^{\frac{n}{2}}}) P_{n}^{-1} := \phi_{n}(a) \in \calC^{2^{\frac{n}{2}} \times 2^{\frac{n}{2}}} =
\calC (2^{\frac{n}{2}})
\label{eq:eq30}
\end{equation}
holds for all $a \in \calC_n.$ Now let $a \in \calC_{n+1}= \calC \{e_1, \ldots, e_{n+1} \}.$ Then $a$ can factor as 
\begin{equation} 
a = a_0 + a_1e_{[n+1]} =  a_0 + e_{[n+1]}a_1, 
\label{eq:eq31}
\end{equation}
where 
$$
a_0, a_1 \in \calC_{n} = \calC\{e_1,\ldots,e_n \}, \quad 
e_{[n+1]}^2 = (-1)^{\frac12(n+1)(n+2)} := r.
$$ 
Moreover define $\oa = a_0 - a_1e_{[n+1]}.$ In that case, 
$D_a = \diag (a I_{2^{\frac{n}{2}}}, \oa I_{2^{\frac{n}{2}}})$ satisfies the following universal
similarity factorization equality
\begin{align}
P_{n+1}D_a P_{n+1}^{-1}
&= \left[ \begin{array}{cc} 
          \phi_{n}(a_0) + \sqrt{r}\phi_{n}(a_1) & 0 \\ 
          0 & \phi_{n}(a_0) - \sqrt{r}\phi_{n}(a_1) 
         \end{array}
  \right] \notag \\[1ex]
&:= \phi_{n+1}(a) \in \, ^2 \calC(2^{\frac{n}{2}}),
\label{eq:eq32}
\end{align}
where
\begin{equation}
P_{n+1} = \frac12 \left[ \begin{array}{cc}
          (1 + \frac{1}{\sqrt{r}}e_{[n+1]})P_{n} & - (\sqrt{r} - e_{[n+1]})P_{n} \\
          \frac{1}{r}(\sqrt{r} - e_{[n+1]})P_{n} & (1 + \frac{1}{\sqrt{r}}e_{[n+1]})P_{n}
                         \end{array} 
                  \right],
\label{eq:eq33}
\end{equation}
and
\begin{equation}
P_{n+1}^{-1} = \frac12 \left[ \begin{array}{cc}
   P_{n}^{-1} (1 + \frac{1}{\sqrt{r}}e_{[n+1]}) & P_{n}^{-1}(\sqrt{r} - e_{[n+1]}) \\  
  -P_{n}^{-1} \frac{1}{r}(\sqrt{r} - e_{[n+1]}) & P_{n}^{-1}(1 +\frac{1}{\sqrt{r}}e_{[n+1]})
                              \end{array} 
                       \right]. 
\label{eq:eq34}
\end{equation}

\medskip
\noindent {\bf Proof.} Note that the commutative rule $be_{[n+1]} = e_{[n+1]}b$ holds for all 
$b \in \calC_{n} = \calC\{e_1, \ldots, e_n \}.$  By applying (\ref{eq:eq30}) to (\ref{eq:eq31}), we obtain
\begin{align*}
P_{n}( aI_{2^{\frac{n}{2}}}) P_{n}^{-1}  &= 
P_{n}( aI_{2^{\frac{n}{2}}}) P_{n}^{-1} + P_{n}(a_1I_{2^{\frac{n}{2}}})P_{n}^{-1}e_{[n+1]} \\
&= \phi_{n}(a_0) + \phi_{n}(a_1)e_{[n+1]} := \psi(a).
\end{align*} 
and 
$$
P_{n}(\oa I_{2^{\frac{n}{2}}}) P_{n}^{-1} =  \phi_{n}(a_0) - \phi_{n}(a_1)e_{[n+1]} 
:= \psi(\oa).
$$ 
Next, setting 
$$ 
V = \frac12 \left[ \begin{array}{cc} 
(1 + \frac{1}{\sqrt{r}} e_{[n+1]})I_{2^{\frac{n}{2}}} & 
    -(\sqrt{r} - e_{[n+1]}) I_{2^{\frac{n}{2}}} \\
 \frac{1}{\sqrt{r}} (\sqrt{r} - e_{[n+1]}) I_{2^{\frac{n}{2}}} &  
     (1 + \frac{1}{\sqrt{r}} e_{[n+1]}) I_{2^{\frac{n}{2}}}
                   \end{array} 
            \right],
$$
$$
V^{-1} = \frac12 \left[ \begin{array}{cc} 
     (1 + \frac{1}{\sqrt{r}} e_{[n+1]}) I_{2^{\frac{n}{2}}} &  
          (\sqrt{r} - e_{[n+1]}) I_{2^{\frac{n}{2}}} \\
     -\frac{1}{r}(\sqrt{r} - e_{[n+1]}) I_{2^{\frac{n}{2}}} &
          (1 + \frac{1}{\sqrt{r}} e_{[n+1]}) I_{2^{\frac{n}{2}}} 
                        \end{array}
                \right],
$$ 
and applying them to $D_a = \diag (\psi(a),\psi(\oa))$ we get 
$$ 
V \left[ \begin{array}{cc}  
          \psi(a) & 0 \\
          0 & \psi(\oa) 
         \end{array} 
  \right] 
=
\left[ \begin{array}{cc} 
        \phi_{n}(a_0) + \sqrt{r} \phi_{n}(a_1) & 0 \\ 
        0 & \phi_{n}(a_0) - \sqrt{r} \phi_{n}(a_1) 
       \end{array} 
\right].
$$
Finally, substituting $\psi(a) = P_{n}( aI_{2^{\frac{n}{2}}}) P_{n}^{-1}$ and
$\psi(\oa) = P_{n}(\oa I_{2^{\frac{n}{2}}})P_{n}^{-1}$ into its left-hand side yields the desired results. \quad $\square$ 

\medskip
\noindent {\bf Theorem 15.}  Suppose that there is an independent invertible matrix $P_{n}$ over $\calC_n = \calC
\{e_1,\ldots,e_{n}\}$ with $n$ even such that
\begin{equation}
 P_{n}( aI_{2^{\frac{n}{2}}}) P_{n}^{-1} := \phi_{n}(a) \in \calC (2^{\frac{n}{2}}) 
\label{eq:eq35}
\end{equation}
holds for all $a \in \calC_n.$ Now, let $a \in \calC_{n+2}=\calC \{e_1,\ldots,e_{n+2}\}.$  Then it can factor as
\begin{equation} 
a = a_0 + a_1e_{[n]}e_{n+1} + a_2 e_{[n]}e_{n+2} + a_3\mu_{n+2},
\label{eq:eq36}
\end{equation}
where 
\begin{gather*}
a_0, a_1, a_2, a_3 \in \calC_n = \calC \{e_1,\ldots,e_n\},\\
(e_{[n]}e_{n+1})^2 = (e_{[n]}e_{n+2})^2 = (-1)^{\frac12 (n+1)(n+2)} = r,\\
\mu_{n+2} = (e_{[n]}e_{n+1})(e_{[n]}e_{n+2}) = -(e_{[n]}e_{n+2})(e_{[n]}e_{n+1}).
\end{gather*} 
In that case, $aI_{2^{\frac{n+1}{2}}}$ satisfies the following universal similarity factorization equality
\begin{align*}
P_{n+2}(aI_{2^{\frac{n+1}{2}}}) P_{n+2}^{-1} & = 
  \left[ \begin{array}{cc}
\phi_{n}(a_0) + \sqrt{r}\phi_{n}(a_1) & r[\phi_{n}(a_2) + \sqrt{r} \phi_{n}(a_3)]  \\             \phi_{n}(a_2) -
\sqrt{r}\phi_{n}(a_3) & \phi_{n}(a_0) - \sqrt{r} \phi_{n}(a_1) 
         \end{array} 
  \right] \\
 & := \phi_{n+2}(a) \in \calC({2^{\frac{n+1}{2}}}), 
\end{align*}
where
$$
P_{n+2} = \frac12 \left[ \begin{array}{cc} 
        (1 + \frac{1}{\sqrt{r}}e_{[n+1]}) P_{n} & 
            (e_{[n]}e_{n+2} - \frac{1}{\sqrt{r}} \mu_{n+2}) P_{n} \\  
        \frac{1}{r}(e_{[n]}e_{n+2} + \frac{1}{\sqrt{r}}\mu_{n+2}) P_{n} & 
            (1 - \frac{1}{\sqrt{r}} e_{[n+1]} ) P_{n} 
                         \end{array}
                  \right],
$$ 
$$
P_{n+2}^{-1} = \frac12 \left[ \begin{array}{cc} 
    P_{n}^{-1}(1 + \frac{1}{\sqrt{r}} e_{[n+1]}) & 
       P_{n}^{-1}(e_{[n]}e_{n+2} - \frac{1}{\sqrt{r}}\mu_{n+2}) \\  
    P_{n}^{-1} \frac{1}{r}(e_{[n]}e_{n+2} + \frac{1}{\sqrt{r}} \mu_{n+2}) &  
       P_{n}^{-1}(1 - \frac{1}{\sqrt{r}} e_{[n+1]}) 
                              \end{array} 
                      \right].
$$
The proof of this result is analogous to that of Theorem 13, and is therefore omitted here. 

\section{Conclusions} 
In this article we have established a set of universal similarity factorization equalities for elements
over the complex Clifford algebra $\Cn.$ These equalities reveal two basic facts about $\Cn:$   
\begin{itemize}
\item[(i)] Each element $a$ in $\Cn$ has a complex matrix representation $\phi_n(a).$ Moreover, a diagonal matrix
constructed by the element $a$ is uniformly similar to its complex matrix representation $\phi_n(a).$
\index{matrix representation!complex}%
\item[(ii)] Conversely, each element $a$ in $\Cn$ could be regarded as an eigenvalue of its complex representation matrix
$\phi_n(a).$ In other words, all complex matrices with the form  $\phi_n(a)$ can uniformly be diagonalized over $\Cn.$  
\end{itemize}

Based on the above two facts, one easily see that almost all known results in complex matrix theory can be extended to
complex Clifford algebras. On the other hand, some problems related to complex matrices can also transform to the
problems related to  Clifford numbers. One such a problem (see \cite{Ab}) is concerning the exponential $e^A$ of a
complex matrix $ A.$ In fact, we can see from Theorems 14 and 15 that for any $a \in \Cn,$ there is 
$$
P_n(e^aI)P_n^{-1} = e^{\phi_n(a)}, 
$$
or 
$$ 
P_n \diag (e^aI, e^{\oa} I ) P_n^{-1} = e^{\phi_n(a)}, 
$$
which implies that the exponential of a complex matrix $\phi_n(a)$ can be determined by the exponential of its
corresponding Clifford numbers.

\end{document}